%% file: HurConc_v2.tex
\documentclass[12pt]{amsart}
\usepackage{amsmath, amssymb, latexsym, amsthm, mathrsfs, color, mathtools, soul, tikz-cd, float}
\usepackage[hyphens]{url}
\usepackage[T1]{fontenc}
\usepackage{hyperref}
\usepackage{lmodern}

\usepackage[top=2in, bottom=1.5in, left=1in, right=1in]{geometry}

\input{shortcuts}

\numberwithin{equation}{thm}

\title{Concentrated sets and the Hurewicz property}

\thanks{The research of the first and the third authors
was funded in whole by the Austrian Science Fund (FWF) [10.55776/I5930 and 10.55776/PAT5730424].
The research of the second author was funded by the National Science Center, Poland Weave-UNISONO call in the Weave programme
Project: Set-theoretic aspects of topological selections 2021/03/Y/ST1/00122
}

\author[V. Haberl]{Valentin Haberl}
\address{Institut f\"ur Diskrete Mathematik und Geometrie, Technische Universit\"at Wien, Wiedner Hauptstrasse 8-10/104, 1040 Wien, Austria.}
\email{valentin.haberl.math@gmail.com}
\urladdr{https://www.tuwien.at/mg/valentin-haberl/}

\author[P. Szewczak]{Piotr Szewczak}
\address{Institute of Mathematics, University of Warsaw, Banacha 2,
02–097 Warsaw, Poland}
\email{p.szewczak@wp.pl}
\urladdr{http://piotrszewczak.pl}

\author[L. Zdomskyy]{Lyubomyr Zdomskyy}
\address{Institut f\"ur Diskrete Mathematik und Geometrie, Technische Universit\"at Wien, Wiedner Hauptstrasse 8-10/104, 1040 Wien, Austria.}
\email{lzdomsky@gmail.com}
\urladdr{https://dmg.tuwien.ac.at/zdomskyy/}



\begin{document}

\maketitle

\begin{abstract}
A set of reals $X$ is $\fb$-concentrated if it has cardinality at least $\fb$ and it contains a countable set $D\sub X$ such that each closed subset of $X$ disjoint from $D$ has size smaller than $\fb$.
We present ZFC results about structures of $\fb$-concentrated sets with the Hurewicz covering property using semifilters.
Then we show that assuming that the semifilter trichotomy holds each $\fb$-concentrated set is Hurewicz and even productively Hurewicz.
We also show that the appearance of Hurewicz $\fb$-concentrated sets under the semifilter trichotomy is somewhat specific and the situation in the Laver model for the consistency of the Borel Conjecture is different.
\end{abstract}

\section{Introduction}

\subsection{Sets of reals and covering properties}
We work in the realm of \emph{sets of reals}, i.e., infinite topological spaces which are homeomorphic to subspaces of the Cantor cube $\Cantor$.
A set of reals $X$ is \emph{Menger}~\cite{Menger24} if for each sequence $\eseq{\cU}$ of open covers of $X$, there are finite sets $\cF_0\sub\cU_0,\cF_1\sub\cU_1,\dotsc$ such that the family $\Un_{n\in\w}\cF_n$ covers $X$.
If in addition the above sets $\cF_n$ are singletons, then the set $X$ is \emph{Rothberger}~\cite{Rot38} and if for each $x\in X$, the sets $
\sset{n}{x\in\Un\cF_n}$ are cofinite then the set $X$ is \emph{Hurewicz}~\cite{Hure25}.
They are central properties considered in the topological selections theory and we have the following implications between them.
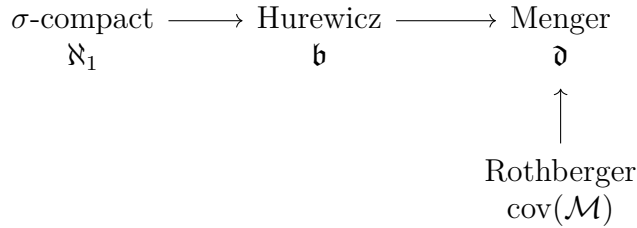
\begin{figure}[H]
\begin{tikzcd}[ampersand replacement=\&,column sep=1cm]
\begin{matrix}\sigma\text{-compact}\\ \aleph_1\end{matrix}\arrow[yshift=2.5mm]{r}\& \begin{matrix}\text{Hurewicz}\\ \fb \end{matrix}\arrow[yshift=2.5mm]{r}\&\begin{matrix}\text{Menger}\\ \fd\end{matrix}\\
{}\& {}\& \begin{matrix}\text{Rothberger}\\ \cov(\cM)\end{matrix} \arrow{u}
\end{tikzcd}
\end{figure}

Each of these properties has its own combinatorial characterization using a structure of the Baire space $\NN$.
For functions $a,b\in\NN$, we write $a\les b$ if the set $\sset{n}{a(n)\leq b(n)}$ is cofinite.
Let $A\sub\NN$.
The set $A$ is \emph{bounded}, if there is a function $b\in\roth$ such that $a\les b$ for all $a\in A$.
A subset of $\NN$ is \emph{unbounded} if it is not bounded.
The set $A$ is \emph{dominating} if for each function $x\in\NN$, there is a function $a\in A$ such that $x\les a$.
The set $A$ is \emph{guessable} if there is a function $g\in\NN$ such that for each function $x\in A$, the sets $\sset{n}{x(n)=g(n)}$ are infinite.
We have the following characterizations of sets of reals with the above properties.

\bthm[{Rec{\l}aw~\cite[Propositions~1, 2, 4]{reclaw}}]\label{thm:rec}
A set of reals $X$ is Hurewicz (Menger, Rothberger) if and only if each continuous image of $X$ into $\NN$ is bounded (not dominating, guessable).
\ethm

Let $\fb$ be the minimal cardinality of an unbounded subset of $\NN$ and $\fd$ be the minimal cardinality of a subset of $\NN$ which is dominating.
By the result of Keremedis \cite[Theorem~2.4.5]{BarJud95}, the minimal cardinality of a subset of $\NN$ which is not guessable is equal to $\cov(\cM)$ which is the minimal cardinality of a family of meager sets in $\NN$ which covers $\NN$.
Theorem~\ref{thm:rec} leads to the conclusion that each of the considered properties has a \emph{critical cardinality} that is the minimal cardinality of a subset of $\NN$ which does not have the given property.
Those critical cardinalities were pointed out in the diagram above.

By the result of Bartoszy\'{n}ski--Shelah~\cite{BaSh01} there is a uniform construction of a Hurewicz set which does not contain a homeomorphic copy of the Cantor cube $\Cantor$.
By the result of Tsaban--Zdomskyy there is also a uniform construction of a Menger set which is not Hurewicz~\cite{sfh}.
In both cases, constructed sets are nontrivial with respect to the given property, i.e., they have cardinalities greater or equal to the critical cardinalities for the given properties, in these cases $\fb$ and $\fd$, respectively.
For the Rothberger property, the situation is more subtle and the existence of a nontrivial Rothberger set is independent of ZFC. 
In the Laver model for the consistency of the Borel conjecture~\cite{Laver} all Rothberger sets are countable.
On the other hand,  if \CH{} holds, an uncountable Rothberger set exists~\cite{Hure25}.

\subsection{Realm of concentrated sets}

In the paper we focus on sets with the following structure.

\bdfn
Let $\kappa$ be an infinite cardinal number.
A set of reals $X$  is \emph{$\kappa$-concentrated on a set} $D\sub X$, if $\card{X}\geq\kappa$ and for each open set $U$ containing $D$, we have $|X\sm U|<\kappa$.
A set is \emph{$\kappa$-concentrated} if it is concentrated on its countable subset.
\edfn

These kinds of sets play a crucial role in considering combinatorial covering properties.
In ZFC there is a $\fd$-concentrated set~\cite[Lemma~17]{Ideals} and each such a set is Menger~\cite[Corollary~1.14]{MHP}.
Also each $\cov(\cM)$-concentrated set is Rothberger but the existence of such a set is independent of ZFC.
In the case of the Hurewicz property it is independent of ZFC, whether each $\fb$-concentrated set is Hurewicz; 
on the one hand, assuming that $\fb=\fd$, there is a $\fb$-concentrated set which is not Hurewicz~\cite{sfh}.
On the other hand, in the Miller model, each $\fb$-concentrated set is Hurewicz (\cite[Theorem~15]{Sch96}, \cite[Theorem~3.1]{semtrich}).

Let $\PN$ be the power set of the set of natural numbers $\w$.
We identify each element of $\PN$ with its characteristic function, an element of $\Cantor$ and in that way introduce a topology in $\PN$.
Let $\PN=\roth\cup\Fin$, where $\roth$ is the family of all infinite subsets of $\w$ and $\Fin$ is the family of all finite subsets of $\w$.
Identifying each element of $\roth$ with increasing enumeration of its elements, we can view $\roth$ as a subset of the Baire space $\NN$.
The topologies in $\roth$ inherited from $\PN$ and $\NN$ are identical.
We refer to the elements of $\roth$ as sets or functions, depending on the contexts.
For practical reasons we mainly work with subspaces of $\PN$.

An example of a Hurewicz $\fb$-concentrated set is a set based on a so-called $\fb$-scale:
a set $X=\sset{x_\alpha}{\alpha<\fb}\sub\roth$ is a \emph{$\fb$-scale} if it is unbounded and for all $\alpha<\beta<\fb$, we have $x_\alpha\les x_\beta$.
A $\fb$-scale exists in ZFC.
If $X\sub\roth$ is a $\fb$-scale, then the set $X\cup\Fin$ is $\fb$-concentrated on $\Fin$ and it is Hurewicz~\cite{BaSh01}, which is a set constructed by Bartoszy\'{n}ski--Shelah, mentioned above.
In fact such a set is even \emph{productively Hurewicz}~\cite[Theorem~6.5]{SPMProd}, i.e., for each Hurewicz set $Y$, the product space $(X\cup\Fin)\x Y$ is Hurewicz.

\subsection{Semifilters}

Considering $\fb$-concentrated sets we need the following notion of semifilters.
For sets $a,b$, we write $a\as b$ if the set $a\sm b$ is finite.
A nonempty set $S\sub\roth$ is a \emph{semifilter}~\cite{BanZdo06} if for any sets $a\in S$ and $b\in\roth$ such that $a\as b$, we have $b\in S$.
Examples of semifilters are the Fr\'{e}chet filter $\cFin$ of all cofinite subsets of $\w$, any free ultrafilter on $\w$ or the full semifilter $\roth$.
Let $S\sub\roth$ be a semifilter.
For functions $a,b\in\roth$, we write $a\leq_S b$ if $\sset{n}{a(n)\leq b(n)}\in S$.
If $S$ is not closed under finite intersections of its elements, then the relation $\leq_S$ need not be transitive; this is the case for $S=\roth$.
The relation $\les$ is the same as $\leq_{\cFin}$ and we write $\leinf$ instead of $\leq_{\roth}$.
A set $A\sub \roth$ is \emph{$\leq_S$-bounded} if there is a function $b\in\roth$ such that $a\leq_S b$ for all functions $a\in A$. 
A set $A\sub\roth$ is \emph{$\leq_S$-unbounded} if it is not $\leq_S$-bounded.
Using the above terminology a set $A\sub\roth$ is $\leq_\cFin$-unbounded if it is unbounded and it is $\leq_\roth$-unbounded if it is dominating.
Let $\bof(S)$ be the minimal cardinality of a $\leq_S$-unbounded set in $\roth$.
We have $\bof(\cFin)=\fb$ and $\bof(\roth)=\fd$.
Moreover $\fb\leq\bof(S)\leq\fd$.

By a recent result of Pawlikowski, Szewczak and Zdomskyy~\cite[Theorem~6.1]{SGom}, in the Miller model any $\fd$-concentrated set has a specific combinatorial structure.
Let $U\sub\roth$ be an ultrafilter.
A set $X\sub\roth$ is a \emph{$U$-scale}, if $|X|\geq\bof(U)$ and for each function $b\in \roth$ there is a set $A\sub X$ of size smaller than $\bof(U)$ such that $b\leq_U x$ for all functions $x\in X\sm A$.
Let $X\sub\roth$.
In the Miller model, the set $X\cup\Fin$ is $\fd$-concentrated on $\Fin$ if and only if there is an ultrafilter $U\sub\roth$ such that $X$ is a $U$-scale.
In fact, in the Miller model, the above set $X$ is a $U$-scale for any ultrafilter $U\sub\roth$.
By the result of Szewczak--Tsaban--Zdomskyy~\cite[Subsection~2.2]{STZ}, this situation from the Miller model is not the case if we assume regularity of $\fd$ and $\fd\leq\fr$, where $\fr$ is the reaping number.
We refer to the work of Blass~\cite{BlassHBK} about cardinal characteristics of the continuum and relations between them.
\mbox{}\\

In Section~2 we present ZFC results about structures of $\fb$-concentrated Hurewicz sets using semifilters.
In Section~3, we show that assuming that the semifilter trichotomy holds\footnote{A statement independent of ZFC, described in detail in Section~3.}, then each $\fb$-concentrated set is Hurewicz and even productively Hurewicz.
In Section~4, we consider $\fb$-concentrated sets in the Laver model for the consistency of the Borel conjecture.
We also show that the appearance of Hurewicz $\fb$-concentrated sets under the semifilter trichotomy is somewhat specific and the situation in the Laver model is different.
 We finish the paper with comments and open problems.

\section{ZFC results}

When we say that a semifilter $S\sub\roth$ is \emph{meager} we mean that it is a meager subset of $\roth$.

\bdfn
A set $X\sub\roth$ is \emph{meager-unbounded} if $\card{X}\geq \fb$ and for each function $b\in\roth$, there are a meager semifilter $S(b)$ and a set $A\sub X$ of size smaller than $\fb$ such that $b\leq_{S(b)} x$ for all functions $x\in X\sm A$.
\edfn

\bthm\label{thm:meagerunbdd}
Let $X\sub\roth$ be a set with $\card{X}\geq\fb$.
Then the following assertions are equivalent.
\be
\item The set $X$ is meager-unbounded.
\item For each $G_\delta$-set $G\sub\PN$ containing $\Fin$, there is a $\sigma$-compact set $F\sub\PN$ such that $F\sub G$ and $\card{X\sm F}<\fb$.
\item The set $X\cup\Fin$ is $\fb$-concentrated on $\Fin$ and for each set $A\sub X$ of size smaller than~$\fb$, the set $(X\sm A)\cup\Fin$ is Hurewicz.
\ee
\ethm

In order to prove Theorem~\ref{thm:meagerunbdd} we need the following notions and auxiliary results.

For natural numbers $n,m$ with $n<m$, let $[n,m):=\sset{i}{n\leq i<m}$.
For sets $x\in\roth$, $S\sub\roth$ and a function $h\in\roth$, define
\[
x/h:=\sset{n\in\w}{x\cap[h(n),h(n+1))\neq\emptyset}\quad\text{and}\quad S/h:=\sset{x/h}{x\in S}.
\]

\blem[{Talagrand~\cite[Theorem~21]{Talagrand}}]\label{lem:tal}
A semifilter $S\sub\roth$ is meager if and only if there is a function $h\in\roth$ such that $S/h$ is the Fr\'echet filter $\cFin$.
\elem

Let $\kappa$ be an uncountable cardinal number and $X\sub\roth$.
The set $X$ is \emph{$\kappa$-unbounded} if $\card{X}\geq\kappa$ and for each function $b\in\roth$ there is a set $A\sub X$ of size smaller than $\kappa$ such that $b\leinf x$ for all functions $x\in X\sm A$.

The next lemma is straightforward. 

\blem\label{lem:unbdd}
Let $\kappa$ be a cardinal with $\mathit{cof}(\kappa)>\w$ and $X\sub\roth$.
Then the set $X\cup\Fin$ is $\kappa$-concentrated on $\Fin$ if and only if the set $X$ is $\kappa$-unbounded.
\elem

Let $\w+1$ be the space with the order topology, the one point compactification of $\w$ and consider $(\w+1)^\w$ with the Tychonoff product topology.
Let $Q$ be the family of all functions $f\in (\w+1)^\w$ with the property that there is $n\in\w$ such that $f\restriction n$ is an increasing function of natural numbers and $f(k)=\w$ for all $k\geq n$.
Then a set $\NNbarup:=\roth\cup Q$ is a closed subset of $(\w+1)^\w$ and thus it is compact.
Moreover, there is a natural homeomorphism between $\PN$ and $\NNbarup$ which identifies a set $a\in\Fin$ with a function $a^\frown\la \w,\w,\dotsc\ra$ from $Q$, where we treat $a$ as an increasing sequence of elements from the set $a$.
In that way we can view elements from $\PN$ as those from $\NNbarup$ and for all elements $a,b\in\PN$ consider a set $[a<b]:=\sset{n}{a(n)<b(n)}$.
E.g., for $a\in\roth$ and $b\in\Fin$ the set $[a<b]$ is cofinite.

For a set $Y$ let $\Pof(Y)$ be the power set of $Y$.
Let $X$ and $Y$ be spaces and $\psi\colon X\to \Pof(Y)$ be a map such that for each $x\in X$, the set $\psi(x)$ is a compact subset of $Y$ and  for each $x\in X$ and an open set $V\sub Y$ such that $\psi(x)\sub V$, there is an open set $U\sub X$ such that $x\in U$ and $\Un\sset{\psi(x')}{x'\in U}\sub V$.
Then the map $\psi$ will be called \emph{compact-valued upper semicontinuous} (\emph{cusco}) and the notation $\psi\colon X\Rightarrow Y$ will be used. Define $\psi[X]:=\Un\sset{\psi(x)}{x\in X}$.

\blem[{\cite[Proposition~4.3]{pMGen}}]\label{lem:cusco}
Cusco maps preserve the Hurewicz property, i.e., if $X$ is a Hurewicz space, $Y$ is a space and $\psi\colon X\Rightarrow Y$ is a cusco map then the set $\psi[X]$ is Hurewicz.
\elem

We have the following topological characterization of the Hurewicz property due to Just--Miller--Scheepers--Szeptycki~\cite{coc2}.

\blem[{Just--Miller--Scheepers--Szeptycki~\cite[Theorem~5.7]{coc2}}]\label{lem:Gdelta}
A set $X\sub\PN$ is Hurewicz if and only if for any $G_\delta$-set $G\sub \PN$ containing $X$, there is a $\sigma$-compact set $F\sub\PN$ such that $X\sub F\sub G$.
\elem

For $a,b\in [\w]^\w$ and a relation $R$ on $\w$ we shall denote by
$[a\,R\,b]$ the set $\sset{n\in\w}{a(n)\,R\,b(n)}$.

\bpf[{Proof of Theorem~\ref{thm:meagerunbdd}}]
(1)$\Rightarrow$(2):
The set $\PN\sm G$ is a $\sigma$-compact subset of $\roth$, and thus there is a function $b\in\roth$ such that $y <^* b$ for all $y\in \PN\sm G$.
Since the set $X$ is meager-unbounded there are a set $A \sub X$ of size smaller than $\fb$ and a meager semifilter $S(b)$ such that $b \leq_{S(b)} x$ for all $x \in X\sm A$.
Since $S(b)$ is meager, by Lemma~\ref{lem:tal}, there is a function $h \in \roth$ such that $S(b)/h$ is the Fr\'{e}chet filter $\cFin$.
Define the meager semifilter
\[
S := \smallmedset{s \in [\w]^\w}{\sset{n \in \w}{s \cap [h(n),h(n+1)) \neq \emptyset} \in \mathrm{cFin}}.
\]
 containing $S(b)$.
In particular, we have $b \leq_S x$ for all $x \in X \sm A$. 

Let
\[
F := \sset{x \in \PN}{b \leq_S x}.
\]
We have that $X\sm A\sub F$ and claim that the set $F$ is $\sigma$-compact. 
Fix $x \in \PN$. 
We have $\sset{n \in \w}{b(n) \leq x(n)} \in S$ if and only if there is $n_x \in \w$ such that for every $n \geq n_x$ there is $i \in [h(n),h(n+1))$ with $b(i) \leq x(i)$.
Let $n_x$ be the minimal number with the above property.
Therefore $F=\Un_{n\in\w}F_n$,
where $F_n := \sset{x \in \PN}{n_x = n}$,
and it is enough to see that the sets $F_n$  are closed for all $n\in\w$:
Fix $y \notin F_n$.
Then there is $k \geq n$ such that for every $i \in [h(k),h(k+1))$ we have $y(i) < b(i)$ and hence $[y \restriction h(k+1)+1]$ is an open neighborhood of $y$ disjoint from $F_n$. 
For every $y \in \PN\sm G$ and $x \in F$ we have 
\[
y<^* b\leq_S x,
\]
and thus $y \leinf x$.
It follows that $F \sub G$.
Since $X \sm F \sub A$, we have $|X \sm F| < \fb$.

(2)$\Rightarrow$(3):
It follows immediately that $X \cup \Fin$ is $\mathfrak{b}$-concentrated on $\Fin$. Fix a set $A\sub X$ of size smaller than $\fb$.
We use a characterization from Lemma~\ref{lem:Gdelta}.
Let $G\sub\PN$ be a $G_\delta$-set containing $(X\sm A)\cup\Fin$.
By (2), there is a $\sigma$-compact set $F\sub\PN$ such that $F\sub G$ and $\card{X\sm F}<\fb$, and thus $\card{X\sm(A\cup F)}<\fb$, too.
The set $X\sm(A\cup F)$ is Hurewicz and it is contained in the set $G$, so by Lemma~\ref{lem:Gdelta}, there is a $\sigma$-compact set $F'\sub G$ such that $X\sm(A\cup F)\sub F'$.
Finally we have $(X\sm A)\cup\Fin\sub F\cup F'\cup\Fin\sub G$, and thus the set $(X\sm A)\cup\Fin$ is Hurewicz.

(3)$\Rightarrow$(1):
Fix $b\in\roth$.
The set $X\cup\Fin$ is $\fb$-concentrated on $\Fin$.
By Lemma~\ref{lem:unbdd}, the set $X$ is $\fb$-unbounded, and thus there is a set $A\sub X$ of size smaller than $\fb$ such that $b<^\infty x$ for all $x\in X\sm A$.
Define 
\[
X_A:=(X\sm A)\cup\Fin
\]
and let $\varphi\colon X_A\to\roth$ be a function such that 
\[
\varphi(x):=[b<x]
\]
for $x\in X_A$.
The map $\varphi$ is well defined (see discussion after Lemma~\ref{lem:unbdd}) and it is continuous.
By the assumption, the set $X_A$
is Hurewicz, and thus the set $Y:=\varphi[X_A]$ is Hurewicz, too.
Let $S(b)$ be a semifilter generated by $Y$, i.e., 
\[
S(b)=\sset{a\in\roth}{\exists y \in Y \enskip (y\as a)}.
\]
We show that $S(b)$ is meager.
For each $n\in\w$, the set $\tilde{Y}_n:=\sset{y\cap [n,\w)}{y\in Y}$ is Hurewicz as a continuous image of $Y$.
Then the set $\tilde{Y}=\Un_{n\in\w}\tilde{Y}_n$ is Hurewicz as a countable union of Hurewicz sets.
Let  $\psi\colon \tilde{Y}\Rightarrow \PN$ be a map such that 
\[
\psi(y):=\sset{z\in\PN}{y\sub z}
\]
for all $y\in\tilde{Y}$.
The map $\psi$ is cusco and we have $S(b)=\psi[\tilde{Y}]$.
By Lemma~\ref{lem:cusco}, the set $S(b)$ is a Hurewicz subset of $\roth$.
By Theorem~\ref{thm:rec}, each Hurewicz subset of $\roth$ is bounded, and thus the set $S(b)$ is meager.
Moreover for all $x\in X\sm A$, we have $b\leq_{S(b)} x$.
It follows that $X$ is meager-unbounded.
\epf

A set $X=\sset{x_\alpha}{\alpha<\fb}\sub \roth$ is a \emph{$\fb$-scale} if the set $X$ is unbounded and for all ordinal numbers $\alpha<\beta<\fb$, we have $x_\alpha\les x_\beta$.
A $\fb$-scale exists in ZFC.

\brem
Let $X\sub\roth$ be a $\fb$-scale.
Then for any set $A\sub X$ with $\card{A}<\mathfrak{b}$, the set $X\sm A$ is also a $\fb$-scale.
By the result of Bartoszy\'{n}ski--Shelah~\cite{BaSh01}, the set $(X\sm A)\cup\Fin$ is Hurewicz.
By Theorem~\ref{thm:meagerunbdd}, the set $X$ is meager-unbounded.
\erem

\brem
Property~(2) from Theorem~\ref{thm:meagerunbdd} is motivated by an observation of Pol--Zak\-rzew\-ski  \cite[Remark~4.2]{pz} who noticed that for a $\fb$-scale $X\sub\roth$ the set $X\cup\Fin$ has this property.
Using this, Lemma~\ref{lem:Gdelta} and the fact that any set of size smaller than $\fb$ is Hurewicz, we get that for any $\fb$-scale $X\sub\roth$, the set $X\cup\Fin$ is Hurewicz.
This argument of Pol--Zakrzewski that for a $\fb$-scale $X\sub\roth$ the set $X\cup\Fin$ is Hurewicz, is much simpler than the original one.
It turns out that the property pointed out by Pol--Zakrzewski implies the following one which is equivalent to a property introduced in\footnote{See Subsection~\ref{subsec:dagger} about the original formulation.}~\cite[Definition~2.1]{SP_LMS}.
\erem

\bdfn\label{dfn:dagger}
A set $X\sub\PN$ satisfies property $(\dagger)$ if for each function $\varphi$ which assigns to each countable set $C\sub X$ a $G_\delta$-set $\varphi(C)\sub\PN$ containing $C$, there is a family $\cK$ of compact sets in $\PN$ such that the family $\cK$ has size smaller than $\fb$, it covers $X$ and refines $\sset{\varphi(C)}{C\text{ is a countable subset of }X}$.
\edfn

A set $X$ is \emph{productively Hurewicz} if for any Hurewicz set $Y$, the product space $X\x Y$ is Hurewicz.
We have also the following result which was originally formulated in a slightly different way, but the proof allows to get the formulation presented here. 

\bthm[{Zdomskyy~\cite[Lemma~2.5]{SP_LMS}}]\label{thm:dagger}
Each subset of $\PN$ satisfying $(\dagger)$ is productively Hurewicz.
\ethm

We get the following corollary from Theorems~\ref{thm:meagerunbdd} and~\ref{thm:dagger}.

\bcor\label{cor:pH}
Let $X\sub\roth$ be a meager-unbounded set.
Then the set $X\cup\Fin$ is productively Hurewicz.
\ecor
\bpf
Let $\varphi$ be a function which assigns to any countable set $C\sub X\cup\Fin$ a $G_\delta$-set $\varphi(C)\sub\PN$ containing $C$.
By Theorem~\ref{thm:meagerunbdd}, there is a $\sigma$-compact set $F$ such that $F\sub\varphi(\Fin)$ and $\card{X\sm F}<\fb$.
Let $F=\Un_{n\in\w}F_n$, where each set $F_n$ is compact.
Then the family 
\[
\sset{\{x\}}{x\in (X\sm F)\cup\Fin}\cup\sset{F_n}{n\in\w}
\]
satisfies the properties from Definition~\ref{dfn:dagger}.
By Theorem~\ref{thm:dagger}, the set $X\cup\Fin$ is productively Hurewicz.
\epf

\brem
Corollary~\ref{cor:pH} can be proven also involving combinatorial methods used by Szewczak--Tsaban~\cite[Theorem~5.4]{ST}.
\erem

\section{Applications under the semifilter trichotomy}

\bdfn
The \emph{semifilter trichotomy} is the statement that for any semifilter $S\sub\roth$, there is a function $h\in \roth$ such that $S/h$ is the Fr\'{e}chet filter $\cFin$, an ultrafilter or the full semifilter $\roth$.
\edfn

The semifilter trichotomy is a statement independent of ZFC which holds, e.g.,  in the Miller model~\cite{mill} and it does not hold under \CH{} (see \cite[Theorem~1]{BanBla06} for a much stronger result).
The aim of this section is to prove the following result.

\bthm\label{thm:semHur}
Assume that the semifilter trichotomy holds.
Let $X\sub\roth$.
Then the following statements are equivalent.
\be
\item The set $X\cup\Fin$ is $\fb$-concentrated on $\Fin$.
\item The set $X\cup\Fin$ is $\fb$-concentrated on $\Fin$ and Hurewicz.
\item The set $X\cup\Fin$ is $\fb$-concentrated on $\Fin$ and it is productively Hurewicz.
\item The set $X$ is meager-unbounded.
\ee
\ethm

We shall prove  Theorem~\ref{thm:semHur} by using several  auxiliary results.

\blem \label{lem_almost_large}
Let $X \sub\roth$ and $\kappa$ be a cardinal with $\mathit{cof}(\kappa)>\w$. 
Assume that $X \cup \Fin$ is $\kappa$-concentrated on $\Fin$ and  $\sset{\cU_n}{n \in \w}$ is a sequence of open covers of $X \cup \Fin$.
Then there is a set $A\sub X$ of size smaller than $\kappa$ and for every $n \in \w$ there is $U_n \in \cU_n$ such that for all $x \in (X \sm A) \cup \Fin$ the set $\sset{n \in \w}{x \in U_n}$ is infinite. 
\elem

\begin{proof}
Enumerate the given sequence of open covers as $\sset{\cU_{\la n,m \ra}}{n,m \in \w}$.
Fix $n \in \w$ and pick for every $m \in \w$ a set  $U_{\la n,m \ra} \in \cU_{\la n,m \ra}$ so that $\sset{U_{\la n,m \ra}}{m \in \w}$ covers the set $\Fin$.
Let 
\[
V_n := \bigcup_{m \in \w} U_{\la n,m \ra} 
\quad\text{ and }\quad A := \bigcup_{n \in \w} X \sm V_n.
\]
Since the set $X \cup \Fin$ is $\kappa$-concentrated on $\Fin$
and  $\mathit{cof}(\kappa)>\w$, the set $A$ has size smaller than $\kappa$.
By definition, we already know that for all $y \in \Fin$ holds $y \in U_{\la n,m \ra}$ for infinitely many $\la n,m \ra$.
If $x \in X \sm A$, then 
\[
x \in \bigcap_{n \in \w} V_n = \bigcap_{n \in \w} \bigcup_{m \in \w} U_{\la n,m \ra}.
\]
Thus, if we enumerate the family $\sset{U_{\la n,m \ra}}{n,m \in \w}$ as $\sset{U_n}{n \in \w}$, we have that for all $x \in (X \sm A) \cup \Fin$ the set $\sset{n \in \w}{x \in U_n}$ is infinite.
\end{proof}

Recall that a set $F\sub\NN$ is \emph{guessable} if there exists 
$g\in\w^\w$ such that $[g=f]$ is infinite for all $f\in F$.

\blem\label{lem:guess}
Let $X \sub\roth$ and $\kappa$ be a cardinal with $\mathit{cof}(\kappa)>\w$. 
Assume that $X \cup \Fin$ is $\kappa$-concentrated on $\Fin$ and $\varphi\colon X  \cup \Fin \to \NN$ is a continuous function.
Then there is a set $A\sub X$ of size smaller than $\kappa$ such that the set $\varphi[(X \sm A) \cup \Fin]$ is guessable. 
\elem

\begin{proof}
Fix $n \in \w$. 
Let us consider the open cover of the Baire space 
\[
\mathcal{V}_n := \sset{V^n_m}{m \in \w},\text{ where }V^n_m := \sset{f \in \roth}{f(n) = m}.
\]
The family $\cU_n := \sset{\varphi^{-1}[V^n_m]}{m \in \w}$ is an open cover of $X \cup \Fin$.
By Lemma~\ref{lem_almost_large}, there are a set $A \sub X$ of size smaller than $\kappa$ and a function $g \in \NN$ such that the family
\[
\sset{\varphi^{-1}[V^n_{g(n)}]}{n \in \w}
\]
fulfills that for every $x \in (X \sm A) \cup \Fin$ the set 
\[
\sset{n \in \w}{x \in \varphi^{-1}[V^n_{g(n)}]} = \sset{n \in \w}{ \varphi(x) \in V^n_{g(n)}} = \sset{n \in \w}{\varphi(x)(n) = g(n)}
\]
is infinite.
Thus, the image $\varphi[(X \sm A) \cup \Fin]$ is guessable by $g$.
\end{proof}

\brem\label{rem:guess}
In Lemma~\ref{lem:guess}, the space $\NN$ can be replaced by $\prod_{n\in\w}A_n$ with the product topology, where $A_n$'s are countably infinite discrete spaces.
Note that in a natural way, the notion of guessability can be considered also in the above product $\prod_{n\in\w}A_n$.
\erem

A set $A\sub\roth$ is a \emph{base} for a set $S\sub \roth$ if 
\[
S=\sset{b\in\roth}{a\as b\text{ for some }a\in A}.
\]


\bprp \label{prop_trichotomy}
Let $X \sub\roth$ and $\kappa$ be a cardinal with $\mathit{cof}(\kappa)>\w$. 
Assume that $X \cup \Fin$ is $\kappa$-concentrated on $\Fin$ and $\varphi\colon X  \cup \Fin \to\roth$ is a continuous function. Then there is a set $A\sub X$ of size smaller than $\kappa$ such that $\varphi[(X \sm A) \cup \Fin]$ is not a base for an ultrafilter nor for the full semifilter $\roth$. 
\eprp

\bpf
Let $Y:= \varphi[X \cup \Fin]$ and fix $y \in  Y$.
Let 
\[
y' := \la \{y(0),y(1)\},\{y(2),y(3),y(4),y(5)\},\dotsc\ra \in \Pi_{n \in \w} [\w]^{2^{n+1}},
\]
where in  $\Pi_{n \in \w} [\w]^{2^{n+1}}$ we consider the product topology and the spaces $[\omega]^{2^{n+1}}$ are discrete.
Then $Y' := \sset{y'}{y \in Y}$ is a homeomorphic copy of $Y$ in $\Pi_{n \in \w} [\w]^{2^{n+1}}$.
So in particular, $Y' = \psi[X \cup \Fin]$ for some continuous function $\psi$.

By Lemma~\ref{lem:guess} and Remark~\ref{rem:guess}, there is a set $A
\sub X$ of size smaller than $\kappa$ such that $\psi[(X \sm A) \cup \Fin]$ is guessable by some function $g \in \Pi_{n \in \w} [\w]^{2^{n+1}}$.
By the definition of $\Pi_{n \in \w} [\w]^{2^{n+1}}$, for every $n \in \w$ we have
\[
|g(n) \sm \bigcup_{i <n} g(i)| \geq 2.
\]
There are disjoint sets $a,b \in [\w]^\w$ such that for every $n \in \w$ we have 
\[
a \cap g(n) \neq \emptyset\quad\text{ and }\quad b \cap g(n) \neq \emptyset.
\]
For each $y \in \varphi[(X \sm A) \cup \Fin]$ we have $|a \cap y| = \w$ and $|b \cap y| = \w$.
Therefore, 
\[
a, \w\sm a \notin \sset{x \in \roth}{\exists y \in  \varphi[(X \sm A) \cup \Fin] \enskip (y \as x)},
\]
and thus the above set, which is a semifilter generated by  $\varphi[(X \sm A) \cup \Fin]$,  is not an ultrafilter nor $\roth$. 
\epf

Let $\fu$ be the minimal cardinality of a base for an ultrafilter in $\roth$.

\blem\label{lem:base}
Let $X\sub\roth$.
\be
\item If $X$ is a base for an ultrafilter $U\sub\roth$ and $A\sub X$ is a set of size smaller than $\fu$, then the set $X\sm A$ is also a base for $U$.
\item  If $X$ is a base for the full semifilter $\roth$ and $A\sub X$ is a set of size smaller than $\fc$, then the set $X\sm A$ is also a base for $\roth$.
\ee
\elem

\bpf
(1) Fix a set $u\in U$.
Since the set $U':=\sset{u\cap u'}{u'\in U}$ is an ultrafilter on $u$, the set $\sset{x\in X}{x\as u}$, a base for $U'$, has cardinality at least $\fu$.
Then there is a set $x\in X\sm A$ such that $x\as u$.
It follows that $X\sm A$ is a base for $U$.

(2) Fix a set $u\in \roth$. Then there is a family $\sset{u_\alpha}{\alpha<\fc}$ of infinite subsets of $u$ such that the sets $u_\alpha\cap u_\beta$ are finite for all $\alpha<\beta<\fc$.
Since the set $A$ has size smaller than $\fc$, there are a set $x\in X\sm A$ and $\alpha<\fc$ such that $x\as u_\alpha\sub u$.
It follows that $X\sm A$ is a base for $\roth$.
\epf

For a set $A\sub\roth$ and a function $g\in\roth$, we write $A\les g$ if $a\les g$ for all $a\in A$.

\begin{proof}[{Proof of Theorem~\ref{thm:semHur}}]
In the light of Theorem~\ref{thm:meagerunbdd} and Corollary~\ref{cor:pH} and since subsets of size $\mathfrak{b}$ containing $\Fin$ of a set that is $\mathfrak{b}$-concentrated on $\Fin$ are still $\mathfrak{b}$-concentrated on $\Fin$, it is enough to prove the implication
(1)$\Rightarrow$(2):
Let $X \cup \mathrm{Fin}$ be $\mathfrak{b}$-concentrated on $\Fin$ and let $\varphi\colon X \cup \mathrm{Fin} \to \roth$ be a continuous function.
By the semifilter trichotomy, there is a function $h \in \roth$ such that $\varphi[X \cup \mathrm{Fin}] / h$ is a base for the Fr\'{e}chet filter, an ultrafilter or $\roth$.
The set $\varphi[X \cup \mathrm{Fin}] / h$ is a continuous image of $X \cup \Fin$.
By Proposition~\ref{prop_trichotomy}, there is a set $A\sub X$ of size smaller than $\fb$ such that $\varphi[(X \sm A) \cup \Fin] / h$ is not a base for an ultrafilter nor for the full semifilter $\roth$.
By Lemma~\ref{lem:base}, the set $\varphi[X\cup\Fin]/h$ is a base for the Fr\'{e}chet filter.
Fix $y\in\varphi[X\cup\Fin]$.
Since the set $\sset{n}{y \cap [h(n),h(n+1)) \neq \emptyset}$ is cofinite, there is a set $c \in\Fin$ such that $y \leq h \sm c$.
Take a function $g \in \roth$ such that $\sset{h\sm c}{ c\in\Fin}\les g$.
Then $ \varphi[X \cup \Fin]\les g$.
By Theorem~\ref{thm:rec}, the set $X \cup \Fin$ is Hurewicz.
\end{proof}

A set $X\sub\PN$ has \emph{strong measure zero} if for any sequence $\seq{\epsilon_n}{n\in\w}$ of positive numbers, there is a sequence $\seq{I_n}{n\in\w}$ of subsets of $\PN$ such that $\diam(I_n)<\epsilon_n$ and the family $\sset{I_n}{n
\in\w}$ covers $X$.
Each Rothberger subset of $\PN$ has strong measure zero.

Let $Y$ be a metrizable space and let $\Iso(Y)$ denote the set of all isolated points of $Y$. 
For the Cantor--Bendixson process on $Y$ we use the following notation: $Y^{(0)} := Y$, $Y^{(\alpha+1)} := Y^{(\alpha)} \sm \Iso(Y^{(\alpha)})$ and for limit ordinals $\lambda$, $Y^{(\lambda)} := \bigcap_{\xi < \lambda} Y^{(\xi)}$. 
The space $Y$ is called \emph{scattered} if every
subspace of $Y$ has an isolated point.
Let $Y$ be a scattered space.
The cardinal number $\height(Y) := \min \sset{\alpha}{Y^{(\alpha)} = \emptyset}$ is the \emph{height} of $Y$.
For $y \in Y$ let $\height(y) := \min \sset{\alpha}{y \in Y \setminus Y^{(\alpha + 1)}}$.
We shall need the following folklore fact.

\bthm\label{thm:scr}
Let $X$ be a metrizable space and let $Y \subseteq X$ be a scattered subspace such that $\height(Y)$ is countable.
Then $Y$ is a $G_\delta$-subset of $X$.
In particular, every scattered subspace of a separable metrizable space $X$ is a $G_\delta$-subset of $X$.
\ethm

\bpf
Let $\beta = \height(Y)$.
Since all $Y^{(\alpha)}$ are closed in $Y$, we can find an increasing sequence $\seq{ O_\alpha}{\alpha \leq \beta}$ of open subsets of $X$ with $O_0 := \emptyset$ and $O_\lambda := \bigcup_{\alpha < \lambda} O_\alpha$ for limit ordinals $\lambda \leq \beta$ such that for all $\alpha < \beta$, we have $O_\alpha \cap Y = Y \setminus Y^{(\alpha)}$.
For all $\alpha < \beta$ let $Z_\alpha := Y \cap (O_{\alpha + 1} \setminus O_\alpha)$. We have
\[
Z_\alpha = (Y \cap O_{\alpha+1}) \setminus (Y \cap O_\alpha) = (Y \setminus Y^{(\alpha+1)}) \setminus (Y\setminus Y^{(\alpha)}) = \sset{y \in Y}{\height(y) = \alpha}.
\]
Since $Z_\alpha = \Iso(Y^{(\alpha)})$, the space $Z_\alpha$ is discrete.
Let us observe that
\[
Y = \bigcup_{\alpha < \beta} \Iso(Y^{(\alpha)}) = \bigcup_{\alpha < \beta} Y \cap (O_{\alpha + 1} \setminus O_\alpha).
\]
Hence we can write
\[
X \setminus Y = (X \setminus O_\beta) \cup \bigcup_{\alpha < \beta} ((O_{\alpha+1} \setminus O_\alpha) \setminus Y) = (X \setminus  O_{\beta}) \cup \bigcup_{\alpha < \beta} ((O_{\alpha+1} \setminus O_\alpha) \setminus Z_\alpha).
\]
Since $X$ is metrizable, every discrete subspace is a $G_\delta$-subset and every open subset is an $F_\sigma$-subset.
By the above it follows that $X \setminus Y$ is an $F_\sigma$-subset.
\epf 

\blem \label{lem:reduction to concentrated on Fin}
Let $Y \sub \roth$ be a $\kappa$-concentrated set for a cardinal $\kappa$ with $\mathit{cof}(\kappa)>\w$. Then there are a set $C \in [Y]^{<\kappa}$ such that $Y \setminus C$ is $\kappa$-concentrated and a set $X \sub \roth$ such that $X \cup \Fin$ is $\kappa$-concentrated on $\Fin$ and $Y \setminus C$ is homeomorphic to $X \cup \Fin$.
\elem

\bpf
Let $A \sub Y$ be a countable set such that $Y$ is $\kappa$-concentrated on $A$.
By the Cantor--Bendixson procedure, we can write $A = B \cup S$, where $B \cap S = \emptyset$, $B$ is a crowded closed subspace of $A$ and $S$ is a scattered space.
Let $\overline{B}$  be the closure of $B$ in $\PN$. 
The set $\overline{B} \cup S$ is a  $G_\delta$-subset of $\PN$: $\overline{B}$ is closed 
in $\PN$, every closed subset of $\PN$ is $G_\delta$, and 
$S$ is $G_\delta$ by Theorem~\ref{thm:scr}.
It follows that $\card{Y \sm (\overline{B} \cup S)} < \kappa$. 
Therefore, $\card{Y \setminus \overline{B}} < \kappa$ and $B \neq \emptyset$.
We get that $\overline{B}$ is homeomorphic to $\PN$ and since $\overline{B} \cap S = \emptyset$, the set $Y \cap \overline{B}$ is $\kappa$-concentrated on $B$.
By the countable dense homogeneity of $\PN$, we can find a homeomorphism $\varphi \colon \overline{B} \to \PN$ such that $\varphi[B] = \Fin$.
Let $X := \varphi[(Y \cap \overline{B}) \sm B] \sub \roth$.
Then $X \cup \Fin = \varphi[Y \cap \overline{B}]$ is $\kappa$-concentrated on $\Fin$ and $\card{Y \sm (Y \cap \overline{B})}< \kappa$.
Put $C:=Y \setminus (Y \cap \overline{B})$.
\epf

\bthm
In the Miller model, each $\fb$-concentrated set is productively Rothberger.
\ethm

\bpf
Let $Y\sub\roth$ be a $\fb$-concentrated set.
By Lemma \ref{lem:reduction to concentrated on Fin} and since $\fb = \w_1$, there are a countable set $C \sub Y$ and  $X\sub\roth$ such that $X\cup\Fin$ is $\fb$-concentrated on $\Fin$ and $X\cup\Fin$ is homeomorphic to $Y \setminus C$.
Let $Z\sub\roth$ be a Rothberger set.
In the Miller model the semifilter trichotomy holds, and thus any Rothberger set is Hurewicz (\cite[Theorem~15]{Sch96}, \cite[Theorem~3.1]{semtrich}).
Then by Theorem~\ref{thm:semHur}, the set $(X\cup\Fin)\x Z$ is Hurewicz.
By the results of Scheepers~\cite[Theorem~1, Lemma~3]{Schsmz} the product space of a Hurewicz strong measure zero set and a strong measure zero set is strong measure zero.
It follows that the set $(X\cup\Fin)\x Z$ is strong measure zero.
By the result of Fremlin--Miller~\cite[Theorem~8]{FrMill}, each Hurewicz strong measure zero set is Rothberger, so the set $(X\cup\Fin)\x Z$ is Rothberger. Since $X\cup\Fin$ is homeomorphic to $Y \setminus C$ and countable unions of Rothberger sets are Rothberger, $Y$ is productively Rothberger.
\epf

\section{Applications in the Laver model}

In this section we show that the appearance of the Hurewicz $\fb$-concentrated sets in the Laver model differs from the one under the semifilter trichotomy.
In particular, we shall show in Corollary~\ref{cor4.3} that we can substantially weaken 
$(3)$ of Theorem~\ref{thm:meagerunbdd} in the Laver model so that the equivalence in Theorem~\ref{thm:meagerunbdd}  still remains true. On the other hand, Proposition~\ref{prp:biconc}
shows some limits of such a weakening in models of $\w_1<\mathfrak b=\mathfrak d$,  including the Laver model.

By the result of the third author we have the following characterization of Hurewicz sets in the Laver model.

\bthm[{Zdomskyy~\cite[Proposition~2.6]{SP_LMS}}]\label{thm:LavHur}
In the Laver model, the Hurewicz property is equivalent to $(\dagger)$.
\ethm

\blem \label{lem4.2}
In the Laver model, let $A\sub Y\sub\roth$.
Assume that for each countable set $C\sub A$, the set $Y\sm C$ is Hurewicz.
Then $Y\sm A$ is Hurewicz.
\elem

\bpf
Fix a $G_\delta$-set $G\sub\PN$ containing $Y\sm A$.
For a set $C\sub Y$, let $C_A:=C\sm G$.
Let $C\sub Y$ be a countable set. 
By the assumption the set $Y\sm C_A$ is Hurewicz, and by Lemma~\ref{lem:Gdelta}, there is a $G_\delta$-set $H_C\sub\PN$ such that $Y\cap H_C=C_A$.
Let $\varphi(C):=H_C\cup G$.
By Theorem~\ref{thm:LavHur}, there is a family $\cK$ of compact sets in $\PN$ of size smaller than $\fb$ which covers $Y$ and $\cK$ refines the family
\[
\sset{\varphi(C)}{C\text{ is a countable subset of }Y}.
\]
Fix $K\in\cK$ and a countable set $C\sub Y$ such that $K\sub H_C\cup G$.
Since $(K\cap Y)\sm C_A$ is a closed subset of the Hurewicz set $Y\sm C_A$, the set $(K\cap Y)\sm C_A$ is Hurewicz.
By Lemma~\ref{lem:Gdelta}, there is a $\sigma$-compact set $F_K\sub\PN$ such that
\[
(K\cap Y)\sm C_A\sub F_K\sub G.
\]
For each $y\in Y\sm A$, there is a set $K\in\cK$ and a countable set $C\sub Y$ such that 
\[
y\in (K\cap Y)\sm C_A\sub F_K.
\]
Thus,
\[
Y\sm A\sub \Un\sset{F_K}{K\in\cK}\sub G
\]
The set $\Un\sset{F_K}{K\in\cK}$ is Hurewicz as a union of less than $\fb$ many Hurewicz sets.
By Lemma~\ref{lem:Gdelta}, there is a $\sigma$-compact set $F\sub\PN$ such that
\[
Y\sm A\sub \Un\sset{F_K}{K\in\cK}\sub F\sub G.
\]
By Lemma~\ref{lem:Gdelta}, the set $Y\sm A$ is Hurewicz.
\epf

\bcor \label{cor4.3}
In the Laver model for each set $X\sub\roth$ of size $\fb$ the following assertions are equivalent.
\be
\item The set $X$ is meager-unbounded.
\item For each countable set $C\sub X$, the set $(X\sm C)\cup\Fin$ is Hurewicz.
\ee
\ecor
\bpf

$(1)\Rightarrow(2):$ This is a direct consequence of the equivalence
$(1)\Leftrightarrow(3)$ in Theorem~\ref{thm:meagerunbdd}.

$(2)\Rightarrow(1)$: By Lemma~\ref{lem4.2} we know that $(X\sm A)\cup\Fin$ is Hurewicz for any $A\sub X$, hence by the equivalence $(1)\Leftrightarrow(3)$ in Theorem~\ref{thm:meagerunbdd} we are left with the task to show that $X$ is $\fb$-concentrated on $\Fin$, i.e., that
$|X\cap K|<\fb$ for any compact $K\sub \roth$.
Given such a set $K$, for any  countable set $C\sub X\cap K$ find a $G_\delta$-subset $G_C$ of $\roth$ such that $G_C\cap(X\cup\Fin)=C$. This is possible because $(X\sm C)\cup \Fin$
is Hurewicz.
Since
$X\cap K= (X\cup\Fin)\cap K$ satisfies $(\dagger)$ being Hurewicz,
the cover 
\[
\sset{G_C}{C\text{ is a countable subset of }X\cap K}
\]
of the set $X\cap K$ has a subcover of size smaller than $\fb$, and hence $|X\cap K|<\mathfrak b.$

\epf

Let $\ici$ be the family of all subsets of $\w$ that are infinite and co-infinite and $\cFin$ be the family of all cofinite subsets of $\w$.
Consider in $\PN$ the canonical metric $d$ inherited from $\Cantor$.
Let $x\in\PN$, $A\sub\PN$ be a nonempty set and $\epsilon>0$.
Define 
\[
\dist(x,A):=\inf\sset{d(x,a)}{a\in A}
\]
and 
\[
\B(A,\epsilon):=\sset{x\in \PN}{d(x,a)<\epsilon\text{ for some }a\in A}.
\]

In the proof of Proposition~\ref{prp:bdHur} below we shall use several times without mentioning the following straightforward 
\bfct\label{clm_count_compact}
Let $K\sub\PN$ be a countable compact set, and for every $x\in K$ let $A_x\sub \PN$ be the range of a sequence convergent to $x$.
Then for every $x\in K$ there exists a finite set $L_x\sub A_x$
such that the set 
\[
K\cup\bigcup\sset{A_x\sm L_x}{x\in K}
\]
is compact.
\efct

\bprp\label{prp:bdHur}
Assume that $\fb=\fd$.
Then there is a set $X\sub\ici$ such that the set $X\cup\Fin\cup\cFin$ is Hurewicz and $\fb$-concentrated simultaneously on $\Fin$ and on $\cFin$.
\eprp
\bpf
We start by constructing an increasing sequence $\seq{K_m}{m\in\w}$
of countable compact subsets of $\Fin\cup\cFin$, and an auxiliary sequence $\seq{S_m}{m\in\w}$ such that $S_m\sub\w^m$ for all $m\in\w$. More precisely,
let 
\[
K_0:=\{\emptyset\},\quad K_1:=K_0\cup\sset{[n,\w)}{n\in\w},\quad S_0=\emptyset,\quad S_1:=\sset{\la n_0\ra}{n_0\in\w}.
\]
For each $s_1\in S_1$, there is $k_{s_1}\in\w$ such that the set
\[
K_2:=K_1\cup\sset{[s_1(0),n_1)}{s_1\in S_1, n_1\geq k_{s_1}}
\]
is compact.
Let 
\[
S_2:=\sset{{s_1}^\frown n_1}{s_1\in S_1, n_1\geq k_{s_1}}.
\]
For each $s_2\in S_2$, there is $k_{s_2}\in\w$ such that the set
\[
K_3:=K_2\cup\sset{[s_2(0),s_2(1))\cup[n_2,\w)}{s_2\in S_2, n_2\geq k_{s_2}}
\]
is compact.
Let 
\[
S_3:=\sset{{s_2}^\frown n_2}{s_2\in S_2, n_2\geq k_{s_2}}.
\]
Fix $m\in\w$ and assume that a compact set $K_m\sub\Fin\cup\cFin$ and a set $S_m$ of increasing sequences of natural numbers of length $m$ have been defined as above.

If $m$ is odd, then for each sequence $s\in S_m$, there is $k_{s}\in\w$ such that the set 
\[
K_{m+1}:=K_m\cup\smallmedset{[s(0),s(1))\cup[s(2),s(3))\cup\dotsb\cup[s(m-1),n)}{s\in S_m, n\geq k_{s}}
\]
is compact since $[s(0),s(1))\cup[s(2),s(3))\cup\dotsb\cup[s(m-1),\w) \in K_m$.
If $m$ is even, then for each sequence $s\in S_m$, there is $k_{s}\in\w$ such that the set 
\[
K_{m+1}:=K_m\cup\smallmedset{[s(0),s(1))\cup[s(2),s(3))\cup\dotsb\cup[s(m-2),s(m-1))\cup [n,\w)}{s\in S_m, n\geq k_{s}}
\]
is compact since $[s(0),s(1))\cup[s(2),s(3))\cup\dotsb\cup[s(m-2),s(m-1)) \in K_m$.
In both cases define
\[
S_{m+1}:=\sset{{s}^\frown n}{s\in S_m, n\geq k_{s}}.
\]
Thus,
\[
K_{m+1}:=K_m\cup\smallmedset{[s(0),s(1))\cup[s(2),s(3))\cup\dotsb\cup[s(m-1),s(m))}{s\in S_{m+1}}
\]
if $m$ is odd, and
\[
K_{m+1}:=K_m\cup\smallmedset{[s(0),s(1))\cup[s(2),s(3))\cup\dotsb\cup[s(m-2),s(m-1))\cup [s(m),\w)}{s\in S_{m+1}}
\]
if $m$ is even.

\bclm\label{clm:hKn}
For every function $h\in\roth$ and sequence $\seq{\epsilon_n}{n\in\w}$ of positive real numbers, there is a set $I\in\ici$ such that for a function $x:=\Un_{i\in I}[h(i),h(i+1))$ we have $\dist(x,K_n)<\epsilon_n$ for all $n\in\w$.
\eclm

\bpf
Pick $i_0\in\w$ with $2^{-h(i_0)}<\epsilon_0$.
Fix $n\in\w$ and assume that an increasing sequence $\la i_0,\dotsc, i_n\ra$ of natural numbers such that $s_n:=\la h(i_0),\dotsc, h(i_n)\ra\in S_{n+1}$ and $2^{-h(i_l)}<\epsilon_l$ for all $l\leq n$ has been defined.
Pick $i_{n+1}> k_{s_n},i_n$ such that $2^{-h(i_{n+1})}<\epsilon_{n+1}$.
Then $h(i_{n+1})>k_{s_n}$, and thus ${s_n}^\frown h(i_{n+1})\in S_{n+2}$.
Put $I:=\Un_{m\in\w}[i_{2m},i_{2m+1})$ and then $x=\Un_{m\in\w}[h(i_{2m}),h(i_{2m+1}))$.
Let $x_0:=\emptyset$, $x_1:=[h(i_0),\w)$.
For each $m\in\w\sm\{0\}$ define

\[
x_{m+1}:=
\begin{cases}
[h(i_0),h(i_1))\cup\dotsb\cup [h(i_{m-2}), h(i_{m-1}))\cup [h(i_m),\w),&\text{ if $m$ is even},\\
[h(i_0),h(i_1))\cup\dotsb\cup [h(i_{m-1}), h(i_{m})),&\text{ if $m$ is odd}.
\end{cases}
\]
We have $x_m\in K_m$ and 
\[
\dist(x, K_m)\leq d(x,x_m)=2^{-h(i_m)}<\epsilon_m
\]
for all $m\in\w$.
\epf

\bclm\label{clm:tilde}
Let $d\in \roth$ be a function with $d(0)\neq 0$. 
Let $\tilde{d}\in\roth$ be a function such that 
\[
\tilde{d}(0):=d(0)\quad\text{and}\quad \tilde{d}(n+1):=d(\tilde{d}(n))
\]
for all $n\in\w$.
If $x\in\roth$ and $x\cap [\tilde{d}(n),\tilde{d}(n+1))=\emptyset$ for some $n\in\w$, then $d(\tilde{d}(n))\leq x(\tilde{d}(n))$.\qedhere
\eclm

Since $\fb=\fd$, there is a dominating scale $\sset{d_\alpha}{\alpha<\fb}\sub\roth$.
Assume that $d_\alpha(0)\neq0$ for all $\alpha<\fb$.
Fix $\alpha<\fb$.
For a set $x\in\PN$, let $x\comp:=\w\sm x$.
Applying Claim~\ref{clm:hKn} to the function $\tilde{d}_\alpha$ and to the sequence $\seq{2^{-d_\alpha(n)}}{n\in\w}$, we get a function $x_\alpha\in\roth$ such that
\begin{equation}\label{eq:inf}
\smallmedset{n}{x_\alpha\cap [\tilde{d}_\alpha(n),\tilde{d}_\alpha(n+1))=\emptyset}, \smallmedset{n}{x_\alpha\comp\cap [\tilde{d}_\alpha(n),\tilde{d}_\alpha(n+1))=\emptyset}\in\roth  
\end{equation}
and 
\begin{equation}\label{eq:Kn}
\dist(x_\alpha,K_n)<2^{-d_\alpha(n)}
\end{equation}
for all $n\in\w$.

Define 
\[
X:=\sset{x_\alpha}{\alpha<\fb}.
\]

\bclm
The set $X\cup\Fin\cup\cFin$ is $\fb$-concentrated on $\Fin$ and on $\cFin$.
\eclm

\bpf
Let $b\in\roth$.
Then there is $\alpha<\fb$ such that $b\les d_\alpha$.
Fix $\beta$ with $\alpha\leq \beta<\fb$.
By Claim~\ref{clm:tilde} and \ref{eq:inf}, we have 
\[
b\les d_\alpha\les d_\beta\leinf x_\beta,x_\beta\comp
\]
It follows that the sets $X$ and $\sset{x_\alpha\comp}{\alpha<\fb}$ are $\fb$-unbounded.
By Lemma~\ref{lem:unbdd}, the sets $X\cup\Fin$ and $\sset{x_\alpha\comp}{\alpha<\fb}\cup\Fin$ are $\fb$-concentrated on $\Fin$.
Since $X\cup\cFin$ is a homeomorphic copy of the set $\sset{x_\alpha\comp}{\alpha<\fb}\cup\Fin$, the set $X\cup\cFin$ is $\fb$-concentrated on $\cFin$.
\epf

\bclm
The set $X\cup\Fin\cup\cFin$ is Hurewicz.
\eclm

\bpf
Let $G\sub\PN$ be a $G_\delta$-set containing $X\cup\Fin\cup\cFin$.
Write $G=\bigcap_{n\in\w}G_n$, where the sets $G_n$ are open and the sequence of $G_n$'s is decreasing.
For each $n\in\w$, there is $f(n)\in\w$ such that $\B(K_n,2^{-f(n)})\sub G_n$.
There is $\alpha<\fb$ such that $f\les d_\alpha$.
Fix $\beta$ with $\alpha\leq\beta<\fb$.
Then $f\les d_\beta$, and thus 
\[
x_\beta\in \B(K_n,2^{-d_\beta(n)})\sub \B(K_n,2^{-f(n)})
\]
for all but finitely many $n\in\w$.
It follows that
\[
\sset{x_\beta}{\alpha\leq\beta<\fb}\sub \bigcup_{n\in\w}\bigcap_{i\geq n}\B(K_i,2^{-f(i)})\sub G
\]
and the set $\bigcup_{n\in\w}\bigcap_{i\geq n}\B(K_i,2^{-f(i)})$ is $\sigma$-compact.
The set $\sset{x_\beta}{\beta<\alpha}\cup\Fin\cup\cFin$ is Hurewicz, so there is a $\sigma$-compact set $F\sub\PN$ such that 
\[
\sset{x_\beta}{\beta<\alpha}\cup\Fin\cup\cFin\sub F\sub G.
\]
We have
\[
X\cup\Fin\cup\cFin\sub \bigcup_{n\in\w}\bigcap_{i\geq n}\B(K_i,2^{-f(i)})\cup F\sub G,
\]
and by Lemma~\ref{lem:Gdelta}, the set $X\cup\Fin\cup\cFin$ is Hurewicz.
\epf

This finishes the proof of Proposition~\ref{prp:bdHur}.\qedhere
\epf

In our next proof we shall use the well-known fact that the Cantor set $\PN$ satisfies the   following strong form of homogeneity: For any countable dense subsets $A_0,A_1,B_0,B_1$
 of $\PN$ such that $A_0\cap A_1=B_0\cap B_1=\emptyset$ there exists a
homeomorphism $h\colon\PN\to\PN$ such that $h[A_0]=B_0$ and $h[A_1]=B_1$. This can be proved ``classically'' by a bit more careful version of the back-and-forth argument 
of Cantor, but also can be derived from the proof of \cite[Theorem~2.1]{Med15}  (attributed there to Baldwin--Beaudoin~\cite{BalBea89}), because the case $\kappa=\w$ 
does not need any extra set-theoretic assumptions since $\mathsf{MA}_\w$
holds in ZFC.

\bprp\label{prp:biconc}
Assume that $\w_1<\fb=\fd$.
Then there is a $\fb$-concentrated Hurewicz set $Y$ such that for each countable set $D\sub Y$ such that $Y$ is $\fb$-concentrated on $D$, and for each countable set $D'\sub Y$ such that $D\sub D'$, there is a countable set $A\sub Y\sm D'$ such that $Y\sm A$ is not Hurewicz.
\eprp

\bpf
Let $\sset{D_\alpha}{\alpha<\w_1}$ be a family of pairwise disjoint countable dense sets in $\roth$.
Fix $\alpha<\w_1$.
Let $X\sub\ici$ be the set from Proposition~\ref{prp:bdHur}.

As noticed above, there is an autohomeomorphism $\varphi_\alpha$ of $\PN$ such that $\varphi_\alpha[\Fin]=\Fin$ and $\varphi_\alpha[\cFin]=D_\alpha$.
 
Let $Y_\alpha:=\varphi_\alpha[X]$ and then the set $\Fin\cup Y_\alpha\cup D_\alpha$ is Hurewicz and a $\fb$-concentrated set on $\Fin$ and $D_\alpha$, simultaneously.
Let 
\[
Y:=\Un_{\alpha<\w_1}(\Fin\cup Y_\alpha\cup D_\alpha).
\]
The set $Y$ is Hurewicz as a union of less than $\fb$ many Hurewicz sets.
The set $Y$ is $\fb$-concentrated on $\Fin$ as a union of less than $\fb$ many sets which are $\fb$-concentrated on $\Fin$.

Fix a countable set $D\sub Y$ such that $Y$ is $\fb$-concentrated on $D$ and let $D'\sub Y$ be a countable set such that $D\sub D'$.
Then there is $\alpha<\w_1$ such that $D'\cap D_\alpha=\emptyset$.
Suppose that $Y\sm D_\alpha$ is Hurewicz. 
By Lemma~\ref{lem:Gdelta}, there is a $\sigma$-compact set $F\sub\PN$ such that 
\[
Y\cap F=Y\sm D_\alpha.
\]

It follows that the set  $\Fin\cup Y_\alpha=F\cap(\Fin\cup Y_\alpha\cup D_\alpha)$ is an $F_\sigma$-subset of  $\Fin\cup Y_\alpha\cup D_\alpha$, a contradiction to the fact that $\Fin\cup Y_\alpha\cup D_\alpha$ is $\mathfrak b$-concentrated on $D_\alpha$.
Put $A:=D_\alpha$.
\epf

\section{Comments and open problems}

\subsection{Around property $(\dagger)$}\label{subsec:dagger}

Originally the property from Definition~\ref{dfn:dagger} was defined formally in a slightly different way.
The following observation shows that our modification is equivalent to the original one.

\bobs
Let $X\sub\PN$.
Then the following assertions are equivalent.
\be
\item The set $X$ satisfies $(\dagger)$.
\item For each function $\varphi$ which assigns to any countable set $C\sub X$ a $G_\delta$-set $\varphi(C)\sub\PN$ containing $C$, there is a family $\cQ$ of countable subsets of $X$ and for each $Q\in \cQ$ there is a $\sigma$-compact set $K_Q\sub\varphi(Q)$ such that the family $\cQ$ has size smaller than $\fb$ and the family $\sset{K_Q}{Q\in\cQ}$ covers $X$. 
\ee
\eobs

\bpf
(1)$\Rightarrow$(2):
Let $\varphi$ be a function which assigns to each countable set $C\sub X$ a $G_\delta$-set $\varphi(C)\sub\PN$ containing $C$.
By (1), there is a family $\cK$ of compact sets in $\PN$ such that the family $\cK$ has size smaller than $\fb$, it covers $X$ and refines 
\[
\sset{\varphi(C)}{C\text{ is a countable subset of }X}.
\]
Then there is a family $\cQ$ of countable subsets of $X$ such that the family $\cQ$ has size smaller than $\fb$ and $\cK$ refines $\sset{\varphi(Q)}{Q\in\cQ}$.
Fix $Q\in\cQ$ and let $\cK_Q$ be the family of all sets $K\in\cK$ such that $K\sub\varphi(Q)$.
The set $\Un\cK_Q$ is Hurewicz as a union of less than $\fb$ many compact sets.
Since $\Un\cK_Q\sub \varphi(Q)$, by Lemma~\ref{lem:Gdelta}, there is a $\sigma$-compact set $K_Q\sub\PN$ such that
\[
\Un\cK_Q\sub K_Q\sub\varphi(Q).
\]
Then the families $\cQ$ and $\sset{K_Q}{Q\in\cQ}$ are as required.

(2)$\Rightarrow$(1): Straightforward.
\epf

\subsection{Open problems}

\bprb Is it true that any Hurewicz $\fb$-concentrated set is productively Hurewicz?
What if \CH{} holds?
\eprb

\bprb
Is the semifilter trichotomy equivalent to the statement that each $\fb$-concentrated set is Hurewicz?
\eprb

\section*{Acknowledgments}
We are grateful to the anonymous referee for his/her careful reading of the manuscript and valuable corrections.

\section*{Declarations}

\subsection*{Publishing policy}
We have read and understood the publishing policy, and submit this
manuscript in accordance with this policy.

\subsection*{Competing interests policy}
The authors declare that they have no competing interests as
defined by Springer, or other interests that might be perceived to influence the results and/or
discussion reported in this paper.

\subsection*{Dual publication} The results/data/figures in this manuscript have not been published
elsewhere, nor are they under consideration
by another publisher.

\subsection*{Authorship}
The corresponding author confirms that he has read the journal policies and
submit this manuscript in accordance with those policies.
Third party material. All of the material is owned by the authors and/or no permissions
are required.
Data availability. Not applicable. This manuscript does not report data generation or
analysis.

\subsection*{Research Funding}
The research of the first and the third authors was funded in whole
by the Austrian Science Fund (FWF) [10.55776/I5930 and 10.55776/PAT5730424]. The
research of the second author was funded by the National Science Center, Poland Weave-
UNISONO call in the Weave programme Project: Set-theoretic aspects of topological selec-
tions 2021/03/Y/ST1/00122.

\subsection*{Authors’ contributions}
All authors whose names appear on the submission have made substantial
contributions to the conception or design of the work; or the acquisition, analysis, or
interpretation of data; or the creation of new software used in the work; drafted the work or
revised it critically for important intellectual content; approved the version to be published;
and agree to be accountable for all aspects of the work in ensuring that questions related to
the accuracy or integrity of any part of the work are appropriately investigated and resolved.

\bibliographystyle{abbrvurl}

\bibliography{bibliography}

\end{document}

%% file: shortcuts.tex
\newcommand{\nc}{\newcommand}

\nc{\thusfar}{\par\bigskip\centerline{\my{--- Edited thus far ---}}\par\bigskip}
\nc{\lei}{\le^\oo}
\nc{\card}[1]{\left|#1\right|}
\nc{\medcard}[1]{\biggl|\,#1\,\biggr|}
\nc{\smallcard}[1]{|\,#1\,|}
\nc{\bds}{bidirectional $\roth$-scale}
\nc{\bbN}{\mathbb{N}}
\nc{\beq}{\begin{eqnarray*}}\nc{\eeq}{\end{eqnarray*}}
\nc{\mbq}{\mb{?}}
\nc{\mb}[1]{{\mbox{\textbf{#1}}}}
\nc{\nop}{$\times$}
\nc{\w}{\omega}
\nc{\fbn}{\!\!\fbox{\!\nop\!}\!\!}
\nc{\yup}{\checkmark}
\nc{\forces}{\Vdash}
\nc{\name}[1]{\dot{#1}}
\nc{\tf}{\my{FINISHED THUS FAR}}
\nc{\FU}{Fr\'echet--Urysohn}
\nc{\gs}{$\gamma$~space}
\nc{\Ga}{\Gamma}\nc{\Om}{\Omega}
\nc{\smallbinom}[2]{\begin{psmallmatrix} #1\\ #2 \end{psmallmatrix}}
\nc{\bgamma}{\smallbinom{\Om}{\Ga}}

\nc{\productive}[2]{\bigl(#1,\allowbreak #2\bigr)^\x}
\nc{\Sel}{\mathsf{S}}
\nc{\sset}[2]{\{\,#1 : #2\,\}}
\nc{\smb}[1]{{\!\!\mb{#1}\!\!}}
\nc{\medset}[2]{{\biggl\{\,#1 : #2\,\biggr\}}}
\nc{\smallmedset}[2]{{\bigl\{\,#1 : #2\,\bigr\}}}
\nc{\set}[2]{{\left\{\,#1 : #2\,\right\}}}
\nc{\seq}[2]{{\la\, #1 : #2\,\ra}}
\nc{\eseq}[1]{#1_0, \allowbreak #1_1, \allowbreak\dots} 
\nc{\cube}{(\Cantor)^\bbN}
\nc{\Match}{\op{Match}}
\nc{\concat}[1]{\hat{\phantom{a}}\langle #1\rangle}
\nc{\poset}{\mathbb{P}}
\nc{\fn}[1]{{\op{Fn}(#1\times\w,2)}}
\nc{\linadd}{\op{linadd}}
\nc{\nonprod}{\non^\x}
\nc{\alephes}{{\aleph_0}}
\nc{\my}[1]{\marginpar{\textcolor{red}{***}}\textcolor{red}{#1}}
\nc{\later}[1]{{\color{green} #1}}
\nc{\BTs}[1]{{\color{green} #1 (BT)}}
\nc{\Cp}{\op{C}_\mathrm{p}}
\nc{\Bp}{\op{B}_p}
\nc{\Pa}[8]{\bibitem{#1} {#2}, \emph{#3}, {#4} \textbf{#5} ({#6}), {#7}--{#8}.}
\nc{\tPa}[5]{\bibitem{#1} {#2}, \emph{#3}, {#4}, to appear.}
\nc{\sPa}[4]{\bibitem{#1} {#2}, \emph{#3}, {#4}, submitted.}
\nc{\Bc}[9]{\bibitem{#1} {#2}, \emph{#3}, in: \textbf{#4} (#5), #6 #7, #8--#9.}
\nc{\fD}{\mathfrak{D}}
\nc{\fX}{\mathfrak{X}}
\nc{\Onbd}{\Op_{\mathrm{nbd}}} 
\nc{\Omnb}{\Om_{\mathrm{nbd}}} 
\nc{\od}{\mathfrak{od}}
\nc{\Setting}[7]{\xymatrix@R=4pt@C=7pt{#1\ar@{-}[r]&#2\ar@{-}[r]&#3\\&#4\ar@{-}[u]\\
#5\ar@{-}[uu]\ar@{-}[r] & #6\ar@{-}[u]\ar@{-}[r] & #7\ar@{-}[uu]}}
\nc{\mx}[1]{\begin{matrix}#1\end{matrix}}
\nc{\plim}{p\txt{-}\lim}
\nc{\Bgp}{{\Z^\bbN}}
\nc{\Cgp}{{{\Z_2}^\bbN}}
\nc{\Cite}[1]{\textbf{[#1]}}
\nc{\Next}[1]{{#1^+}}
\nc{\cFin}{\mathrm{cFin}}
\nc{\intvl}[2]{{[#1(#2),\allowbreak #1(#2\!+\!1))}}
\nc{\Bdd}{\mathbf{B}}
\nc{\Dfin}{\mathfrak{D}_\mathrm{fin}}
\nc{\grbl}{{\mbox{\textit{\tiny gp}}}}
\nc{\bbP}{\mathbb{P}}
\nc{\BOfat}{\B_{\Om_{\mathrm{fat}}}}
\nc{\Bgood}{\B_{\mathrm{good}}}
\nc{\compactN}{\cl{\mathbb{N}}}
\nc{\blocks}[2]{\op{cl}_{#2}(#1)}
\nc{\blocksplus}[2]{\op{cl}^+_{#2}(#1)}
\nc{\arx}[1]{\texttt{http://arxiv.org/math/#1}}
\nc{\bq}{\begin{quote}}
\nc{\eq}{\end{quote}}
\nc{\cl}[1]{\overline{#1}}
\nc{\CH}{the Continuum Hypothesis}
\nc{\MA}{Martin's Axiom}
\nc{\Bfat}{\B_\mathrm{fat}}
\nc{\inv}{^{-1}}
\nc{\Cantor}{{2^\w}}
\nc{\bP}{\mathbf{P}}
\nc{\bof}{\op{\fb}}
\nc{\dof}{\op{\fd}}
\nc{\bofF}{\bof(\cF)}
\nc{\sr}[3]{\underset{\mbox{#3}}{\mbox{#1}}}
\nc{\gp}{\binom{\Om}{\Ga}}
\nc{\gpsmall}{\mbox{$\gp$}}
\nc{\gig}{\gimel}
\nc{\gns}{\sone(\Om,\gig)}
\nc{\nsr}[2]{#1}
\nc{\Srg}{{\mathbb{S}}}
\nc{\Srgs}{{\mathbb{S}^*}}
\nc{\NN}{{\w^{\w}}}
\nc{\ZN}{{\Z^{\bbN}}}
\nc{\NNup}{{\bbN^{\uparrow\bbN}}}
\nc{\NNbarup}{{(\w+1)^{\uparrow\w}}}
\nc{\Pof}{\op{P}}
\nc{\PN}{{\Pof(\w)}}
\nc{\roth}{{[\w]^{\w}}}
\nc{\Fin}{\mathrm{Fin}}
\nc{\ici}{[\w]^{ \w, \w}}
\nc{\Inc}{{\compactN^{\uparrow\bbN}}}
\nc{\powInc}[1]{{\big(\Inc\big)^{#1}}}
\nc{\powFin}[1]{{\big(\Fin\big)^{#1}}}
\nc{\powPN}[1]{{\big(\PN\big)^{#1}}}
\nc{\NcompactN}{{\compactN^\bbN}}
\nc{\Uarrow}{\smash{\big\uparrow}}
\nc{\LE}{\preccurlyeq}
\nc{\GE}{\succcurlyeq}
\nc{\op}{\operatorname}
\nc{\im}{\op{im}}
\nc{\Span}{\op{span}}
\nc{\maxfin}{\op{maxfin}}
\nc{\ran}{\op{range}}
\nc{\iso}{\cong}
\nc{\Madd}{{\M}^\star}
\nc{\cI}{\mathcal{I}}
\nc{\cJ}{\mathcal{J}}
\nc{\scrA}{\mathscr{A}}
\nc{\scrB}{\mathscr{B}}
\nc{\scrC}{\mathscr{C}}
\nc{\scrD}{\mathscr{D}}
\nc{\scrF}{\mathscr{F}}
\nc{\scrK}{\mathscr{K}}
\nc{\A}{\forall}
\nc{\B}{\mathrm{B}}
\nc{\cB}{\mathcal{B}}
\nc{\bB}{\mathbf{B}}
\nc{\BS}{\mathbf{B}(\mathcal{S})}
\nc{\BF}{\mathbf{B}(\mathcal{F})}
\nc{\BU}{\mathbf{B}(\mathcal{U})}
\nc{\cSp}{\mathcal{S}^+}
\nc{\cFp}{\mathcal{F}^+}
\nc{\cUp}{\mathcal{U}^+}

\nc{\BG}{\B_\Ga}
\nc{\BL}{\B_\Lambda}
\nc{\BT}{\B_\Tau}
\nc{\BTstar}{\B_{\Tau^*}}
\nc{\BO}{\B_\Om}
\nc{\DO}{\cD_\Om}
\nc{\KO}{\cK_\Om}
\nc{\CG}{C_\Ga}
\nc{\CL}{C_\Lambda}
\nc{\CT}{C_\Tau}
\nc{\CTstar}{C_{\Tau^*}}
\nc{\CO}{C_\Om}
\nc{\COgp}{C_{\Om^{\grbl}}}
\nc{\CLgp}{C_{\Lambda^{\grbl}}}
\nc{\BOgp}{\B_{\Om}^{\grbl}}
\nc{\BLgp}{\B_{\Lambda^{\grbl}}}
\nc{\sfC}{\mathsf{C}}
\nc{\sfD}{\mathsf{D}}
\nc{\bD}{\mathbf{D}}
\nc{\Tau}{\mathrm{T}}
\nc{\cA}{\mathcal{A}}
\nc{\cK}{\mathcal{K}}
\nc{\cD}{\mathcal{D}}
\nc{\cF}{\mathcal{F}}
\nc{\cS}{\mathcal{S}}
\nc{\cT}{\mathcal{T}}
\nc{\cG}{\mathcal{G}}
\nc{\cY}{\mathcal{Y}}
\nc{\J}{\mathcal{J}}
\nc{\cL}{\mathcal{L}}
\nc{\cM}{\mathcal{M}}
\nc{\cN}{\mathcal{N}}
\nc{\cH}{\mathcal{H}}
\nc{\cO}{\mathcal{O}}
\nc{\Op}{\mathrm{O}}
\nc{\rmA}{\mathrm{A}}
\nc{\rmF}{\mathrm{F}}
\nc{\rmB}{\mathrm{B}}
\nc{\rmD}{\mathrm{D}}
\nc{\rmP}{\mathrm{P}}
\nc{\cC}{\mathcal{C}}
\nc{\cP}{\mathcal{P}}
\nc{\bbQ}{\mathbb{Q}}
\nc{\bbR}{\mathbb{R}}
\nc{\cU}{\mathcal{U}}
\nc{\cQ}{\mathcal{Q}}
\nc{\Un}{\bigcup}
\nc{\cV}{\mathcal{V}}
\nc{\cW}{\mathcal{W}}
\nc{\Z}{{\mathbb Z}}
\nc{\Impl}{\Rightarrow}
\long\def\forget#1\forgotten{\marginpar{\textcolor{green}{Forgetting...}}}
\nc{\ft}{\mathfrak{t}}
\nc{\fb}{\mathfrak{b}}
\nc{\fc}{\mathfrak{c}}
\nc{\fd}{\mathfrak{d}}
\nc{\fg}{\mathfrak{g}}
\nc{\oo}{\infty}
\nc{\fr}{\mathfrak{r}}
\nc{\fk}{\mathfrak{k}}
\nc{\bidi}{\mathfrak{bidi}}
\nc{\fu}{\mathfrak{u}}
\nc{\fh}{\mathfrak{h}}
\nc{\fp}{\mathfrak{p}}
\nc{\fj}{\mathfrak{j}}
\nc{\fs}{\mathfrak{s}}
\nc{\x}{\times}
\nc{\Iff}{\Leftrightarrow}
\newcommand\comp{^{\text{\tt c}}}
\nc{\nin}{\notin}
\nc{\cat}{\hat{\ }}
\nc{\sub}{\subseteq}
\nc{\spst}{\supseteq}
\nc{\sm}{\setminus}
\nc{\as}{\subseteq^*}
\nc{\les}{\le^*}
\nc{\leinf}{\le^{\infty}}
\nc{\leS}{\le_S}
\nc{\leF}{\le_{\mathcal{F}}}
\nc{\leU}{\le_{U}}
\nc{\rest}{\restriction}

\nc{\la}{\langle}
\nc{\ra}{\rangle}
\nc{\E}{\exists}
\nc{\dom}{\op{dom}}
\nc{\cov}{\op{cov}}
\nc{\add}{\op{add}}
\nc{\cof}{\op{cof}}
\nc{\cf}{\op{cf}}
\nc{\non}{\op{non}}
\nc{\unif}{\op{non}}
\nc{\COV}{\op{COV}}
\nc{\ADD}{\op{ADD}}
\nc{\COF}{\op{COF}}
\nc{\NON}{\op{NON}}
\nc{\impl}{\to}
\nc{\Lp}{\mathcal{L_\p}}
\nc{\Wlog}{without loss of generality}
\newtheorem{thm}{Theorem}[section]
\nc{\bthm}{\begin{thm}} \nc{\ethm}{\end{thm}}
\newtheorem{prop}[thm]{Proposition}
\nc{\bprp}{\begin{prop}} \nc{\eprp}{\end{prop}}
\newtheorem{fact}[thm]{Fact}
\nc{\bfct}{\begin{fact}} \nc{\efct}{\end{fact}}
\newtheorem{prob}[thm]{Problem}
\nc{\bprb}{\begin{prob}} \nc{\eprb}{\end{prob}}
\newtheorem{lem}[thm]{Lemma}
\nc{\blem}{\begin{lem}} \nc{\elem}{\end{lem}}
\newtheorem{claim}[thm]{Claim}
\nc{\bclm}{\begin{claim}} \nc{\eclm}{\end{claim}}
\newtheorem{cor}[thm]{Corollary}
\nc{\bcor}{\begin{cor}} \nc{\ecor}{\end{cor}}
\newtheorem{conj}[thm]{Conjecture}
\nc{\bcnj}{\begin{conj}} \nc{\ecnj}{\end{conj}}
\theoremstyle{definition}
\newtheorem{defn}[thm]{Definition}
\nc{\bdfn}{\begin{defn}} \nc{\edfn}{\end{defn}}
\newtheorem{obs}[thm]{Observation}
\nc{\bobs}{\begin{obs}} \nc{\eobs}{\end{obs}}
\theoremstyle{remark}
\newtheorem{rem}[thm]{Remark}
\nc{\brem}{\begin{rem}} \nc{\erem}{\end{rem}}
\newtheorem{cnv}[thm]{Convention}
\nc{\bcnv}{\begin{cnv}} \nc{\ecnv}{\end{cnv}}
\newtheorem{exam}[thm]{Example}
\nc{\bexm}{\begin{exam}} \nc{\eexm}{\end{exam}}
\nc{\bpf}{\begin{proof}} \nc{\epf}{\end{proof}}
\nc{\be}{\begin{enumerate}}
\nc{\ee}{\end{enumerate}}
\nc{\bi}{\begin{itemize}}
\nc{\bimy}{\my{\begin{itemize}}
\nc{\eimy}{\end{itemize}}}
\nc{\itm}{\item}
\nc{\ei}{\end{itemize}}
\nc{\Subsection}[1]{\goodbreak\subsection*{#1}}

\nc{\sone}{\mathsf{S}_1}
\nc{\sfin}{\mathsf{S}_\mathrm{fin}}
\nc{\ufin}{\mathsf{U}_\mathrm{fin}}
\nc{\Split}{\mathrm{Split}}

\nc{\gone}{\mathsf{G}_1}   
\nc{\Succ}{\mathrm{S}} 
 \nc{\gfin}{\mathsf{G}_\mathrm{fin}}
\DeclareMathOperator{\dist}{dist}

\DeclareMathOperator{\eexists}{\exists}
\DeclareMathOperator{\fforall}{\forall}
\nc{\Exists}[1]{\bigl(\eexists #1\bigr)}
\nc{\Forall}[1]{\fforall #1\ }
\nc{\Foralm}[1]{\fforall^* #1\ }
\nc{\End}[1]{#1}
\nc{\supp}{\op{supp}}
\nc{\bfP}{\mathbf{P}}
\nc{\diam}{\op{diam}}
\nc{\Iso}{\op{Iso}}
\nc{\height}{\op{height}}